\documentclass{amsart}

\usepackage{amsmath,amsthm,graphicx,hyperref}

\usepackage{randtext}

\pdfoutput=1

\usepackage{epigraph}
\newcommand{\AM}{\textsc{AutoMath}}

\title{Mathematics and Mathematics Education in the 21st Century}

\author{Alexandre Borovik$^{1}$}
\address{University of Manchester, UK}
\email{alexandre$\gg$at$\ll$borovik.net}
\thanks{$^{1}$ ORCID 000-0003-0808-8687}

\author{Zoltan Kocsis}
\address{ University of Stirling, UK}
\email{zak$\gg$at$\ll$cs.stir.ac.uk}

\author{Vladimir Kondratiev$^{2}$}\thanks{$^{2}$ ORCID 0000-0002-9792-2698}
\address{Higher School of Economics, Moscow, Russia}
\email{kondratjew239$\gg$at$\ll$gmail.com}

\thanks{\copyright  2022 by the authors. Submitted
for possible open access publication
under the terms and conditions
of the Creative Commons Attribution (CC BY) license (\href{https://
creativecommons.org/licenses/by/
4.0/}{https://
creativecommons.org/licenses/by/
4.0/}). }

\begin{document}

\maketitle

\begin{abstract}
Mathematics enters the period of change unprecedented in its history, perhaps even a revolution: a switch to use of computers as assistants and checkers in production of proofs. This requires rethinking traditional approaches to mathematics education which is struggling through a crisis of its own, socio-economic and political by its nature. The mathematical community faces Pandora's box of problems, which, surprisingly, are not usually discussed in any connected form. The present paper attempts to address this issue in a bit more joint and cohesive way.
\end{abstract}

\section{Introduction}

\subsection{What this paper is about and for whom it is written}

In this paper, we try to explain that mathematics as a scientific discipline enters the period of crisis. In plain words, research level mathematics is becoming too complex for human comprehension. This has serious implications for mathematics education, and our paper is written for everyone who cares about it.

\subsection{Disclaimers}
Hopefully the following clarifications  will preempt some questions to the authors.
\begin{quote}
\begin{itemize}
 \item The aim of this paper is to pose questions and attract attention to issues in recent (not necessary positive) developments  in mathematics, mathematics education, and the socio-political environments in which they exist.

 \item Talking about mathematics education, we do not focus on any particular country. The authors have experience of studying,  and/or doing research, and/or teaching in Australia, England, Hungary, Russia, Scotland, Turkey,  USA. We feel that what we write here applies to almost any industrial democracy.

 \item This paper is not a survey. A few examples we give serve  only as indicators of possible directions of development.

 \item Mentioning any software product, teaching method, etc.  is not an endorsement.

 \item \textbf{All authors write in their private capacity and their views do not necessarily represent positions of their employers, or any other person, corporation, organisation, or institution.}

\end{itemize}
\end{quote}

\section{The crisis in mathematics}
\label{sec:Crisis}

About 40 years ago it was estimated
that between
100 000 and 200 000 new theorems were published every year in mathematical journals around the world \cite{Davis-Hersh1980} --- this number could only have increased since then. A mathematical theorem, as a rule, explicitly refers to other theorems and definitions and is integrated into the huge system of
mathematical knowledge. This system remains unified, tightly connected,  and cohesive; despite the fantastic diversity of mathematics, it also has almost incomprehensible unity.

Mathematics continues to grow, and if you look around, you see that mathematical results and concepts involved in practical applications are much deeper, more abstract and more difficult than ever before. And we have to accept that the mathematics hardwired and coded, say, in a smartphone, is beyond the reach of the vast majority of graduates from  mathematics departments in many universities in the world.

The cutting edge of mathematical research moves further and further away from the stagnating mathematics education. From the point of view of an aspiring PhD student, mathematics looks like New York in the \v{C}apek Brothers’ book \emph{A Long Cat Tale} \cite[p. 44]{Capek1927}:
\begin{quote}\small
And New York -- well, houses there are so tall that they can’t even finish building them.
Before the bricklayers and tilers climb up them on their ladders, it is noon, so they eat their
lunches and start climbing down again to be in their beds by bedtime. And so it goes on day
after day.
\end{quote}
This is mathematics, as seen by many PhD students -- they have no time to climb to the top. In England, most PhD theses, by necessity, have to solve artificial, learner's level problems which contribute next to nothing to the advancement of mathematics.

But Josef and Karel \v{C}apek also formulated an answer to this conundrum: they coined the word `\emph{robot}' (for Karel's play \emph{R.U.R.}, premiered  2 January 1921) for a specific socio-economic phenomenon: a device or machine whose purpose is to replace a human worker.

Mathematics badly needs its own specialised mathematical robots -- first of all, for checking proofs, which are becoming impossibly long and difficult. One of the more notorious examples of length and difficulty is the Classification of the Finite Simple Groups (CFSG), one of the central results of 20th century algebra. In particular, the CFSG controls the behavior of symmetry in the entirety of finite mathematics---everything that deals with finite objects and structures. It is worth remembering that, no matter how huge and complex  they may be, protein molecules and computer processors still remain finite structures.

The original proof of the CFSG, at that time still with some holes and gaps, was spread over more than 100 journal papers of total length about $15\, 000$ pages. A proper and structured  proof  is being published, volume by volume, since 1994 \cite{CSFG}. To date, 9 volumes out of the originally estimated 12 have been published, and tenth volume is in preparation. There are at most a couple dozen people in the world who can read and understand this proof. The youngest of them is Inna Capdeboscq, one of the authors of the latest volume \cite{CGLS2021}; very soon she will be (or, which is already very likely, is) the only non-retired mathematician who understands the entire proof of the CFSG.

Dozens of similar examples can be given: the length and complexity of many very important proofs are now reaching the limits of human comprehension. We have to admit that mathematics faces an existential crisis.
Indeed, what is the purpose of a  mathematical proof that hardly anyone will ever be willing and able to read and check?

Without pivoting to a systematic use of computer-based proof assistants, and corresponding changes to the way mathematics is taught and the way mathematical research is published, the field will not be able to face its new challenges, including those arising in biology, information technology, or data science. Moreover, it is likely to enter a spiral of decay.

\section{The future role of proof assistants}
\label{sec:Proof-assistants}

Baez~\cite{baez-talk} lists "making mathematics more computer-friendly" and "making mathematics more human-friendly"\footnote{To enable mathematicians to acquire the necessary education that allows them tackle the endogenous and exogenous challenges posed by the 21st century.} as two of his five anticipated changes to Mathematics in the 21st century. He points out that these two items are two sides of the same coin, and that the former should have applications in the service of the latter. "Ultimately, all results should be computer-verified, easy to access world-wide, and annotated in various ways."

Simpson's 1998 counterexample to the Homotopy Hypothesis~\cite{simpson-homotopy-theory} overturned an earlier, up to that point widely accepted argument of Kapranov and Voevodsky~\cite{kapranov-homotopy-hypothesis}. Voevodsky's long struggle to identify the error in his own argument led him to formulate the Univalent Foundations program~\cite{voevodsky-univalent-year}, which ultimately resulted in significant advances in the computer verification (mechanization) of proofs, including the development of an entire field of mathematics now known as Homotopy Type Theory.

While the earliest work on computer-assisted mathematics focused on \textit{automated} theorem provers (software that searches for proofs of assertions fully automatically), in the long run, the development of modern, user-friendly \textit{interactive} theorem proving software (\textit{proof assistants}) was the technological breakthrough that enabled the formalization and computer verification of non-trivial amounts of contemporary mathematical research. The purpose of interactive theorem provers is not the automatic generation of proofs, but the verification of proof scripts written by the user (i.e.~the mathematician) in a formal proof language. Proof assistants alleviate many of the crisis factors identified in Section~\ref{sec:Crisis}:

\subsection{Climbing to the top}

Proof assistants help students ``climb to the top'' in multiple ways. Formalized proof bases inherently constitute \textit{explorable examples}, allowing the learner to play with a proof's assumptions, and get immediate feedback about the consequences. Indeed, proof assistants provide major speed-ups in undergraduate education, by enabling immediate feedback. Incorrect proofs are rejected right away, and the error message points the user to the erroneous/incomplete parts of their arguments. For example, if the student wishes to apply a certain theorem, the proof assistant will make sure that the required hypotheses have indeed been shown. This results in accelerated learning when compared to the usual system of mathematics education, where an incorrect assumption or incomplete argument may only be caught by an educator weeks later, long after the student stopped thinking about the question, and at a time when the course has moved on~\cite{nipkow-semantics}.

Proof assistants provide further advantages at the graduate level. On one hand, as a general rule, graduate students are much more capable of evaluating the correctness of their own work independently than undergraduates. On the other hand, for graduates, ``climbing to the top'' (as explained in Section~\ref{sec:Crisis}) requires attention from domain experts, whose time is in even higher demand. The search tools of proof assistants, combined with the contextual information available in mathematical libraries, empower graduate students to answer many questions\footnote{A whimsical example: is Schultz's proof~\cite{schultz-wiles-mo} of the irrationality of $\sqrt[3]{2}$, by reduction to Wiles' theorem, essentially circular? Here is this proof: If $2^{\frac{1}{3}} = \frac{p}{q}$ then $p^3=q^3 + q^3$, which contardicts Wiles' Theorem.} that would otherwise require expert-level knowledge of the literature.

The process of integrating proof assistants into mathematics educations is already underway: e.g.~undergraduate students at Imperial College London, working with Kevin Buzzard, are learning to use proof assistants and are making contributions to the formalization of increasingly sophisticated mathematical content.

\subsection{The CFSG Scenario}

The adoption of proof assistants enables the creation of unified formalized libraries of mathematics, such as Lean's \texttt{mathlib}~\cite{lean-library}. These formalized libraries, combined with sophisticated search tools, reduce the probability of scenarios similar to the situation with CFSG described in Section~\ref{sec:Crisis}, where the knowledge base behind huge projects is rapidly lost, and cannot be recovered due to the size and the spread-out nature of the literature. In an interconnected, annotated proof base, results that see active use inherently see active maintenance.

\subsection{New values on new tablets}

Last, but not least, systematic use of proof assistants enables work on \textit{kinds of mathematics} whose pursuit would have been impossible, or prohibitively expensive, otherwise.

It has been known since the late 20th century that the correctness of certain results that make heavy use of computer calculations can be exorbitantly difficult to verify without proof assistants. Hales' proof of the Kepler conjecture on dense sphere packings, which relied on a large number of such calculations, provides a notorious example: after its submission to the \textit{Annals of Mathematics} in 1998, the proof spent four years in peer review, which ended with an inconclusive result~\cite{nature-szpiro}. MacPherson, then editor of the Annals, wrote:
\begin{quote}
    ``The news from the referees is bad, from my perspective. They have not been able to certify the correctness of the proof, and will not be able to certify it in the future, because they have run out of energy to devote to the problem.''
\end{quote}

Ultimately, this led Hales to start the \texttt{Flyspeck} project, which successfully produced a computer-verified proof of the conjecture in 2015~\cite{hales-flyspeck}.

Even proofs that were obtained without computer assistance benefit from the use of interactive theorem provers. In late 2020, Scholze~\cite{scholze_2020} wrote a blog post (titled Liquid Tensor Experiment), in which he challenged the proof assistants community to formalize a difficult foundational theorem in Condensed Mathematics, a subject pioneered by Scholze in joint work with Clausen. In less than six months, the resulting formalization effort, led by Johan Commelin and using the Lean proof assistant, verified the entire argument that Scholze was unsure about, and demonstrated that interactive theorem provers have reached a level of development where they can be used to formally verify difficult original research within a reasonable time span~\cite{nature-scholze-2021}. Moreover, the formalization answered an open question from Scholze's original lecture notes, and led to new mathematical developments, among them a construction of a variant Breen-Deligne resolution that avoids the use of stable homotopy theory and, according to Scholze himself, ``may become a standard tool in the coming years''~\cite{scholze_2021}.

Finally, proof assistants enable certain powerful styles of mathematical reasoning that were previously too difficult to check not just because of the necessity of expert atteniton, but because they are unusually susceptibile to accidental mistakes due to reliance on syntactic restrictions on set formation and induction principles that have universal validity in the rest of mathematics. For example, one of the authors \cite{Kocsis2019} observed that many arguments in nonstandard analysis rely on such syntactically restricted set formation principles, and deployed proof assistants to verify novel mathematical results in the subject.

Due to the difficulty and tedium of syntactic checks, before the advent of computer assistance, mathematicians gravitated towards foundational theories that did not require them. Indeed, Karagila~\cite{karagila-zf-over-nfu} identifies avoiding the need for such syntactic verification as a significant advantage that Zermelo-Fraenkel set theory held over its foundational alternatives throughout the 20th century:

\begin{quote}
``You need to keep verifying whether the defining formulas are stratified, and if you can or cannot use them to define new sets. [\dots ] I want to argue something about the collection of topologies of the real numbers which are Hausdorff. I don't need to argue that the defining formula of being a Hausdorff topology over $\mathbb{R}$ has some syntactical properties and show that this collection is a set. I simply state it and move along.''
\end{quote}

Thus, proof assistants enable the advancement of fields of mathematics where syntactic checks cannot be disposed with. This includes the foundational field of Homotopy Type Theory (which relies on similar syntactic checks for e.g. universe levels), whose establishment was intertwined with advances in proof assistant technology, and which is probably the only current field of mathematics where formally verified arguments outnumber purely informal developments.

\subsection{The choice of the name}

The new proof-assistants-based mathematics need a catchy name. After some discussions, we agreed that Automatic Mathematics and abbreviated as \AM\ would be the right choice. Gregory Cherlin kindly offered an advertising jingle:
\begin{quote}
\textit{ Take the path of AutoMath\\
 and cut a swath across your math! }
\end{quote}
The name could be seen as a tribute to Nicolaas Govert de Bruijn (1918--2012):  \textit{AutoMath}\  was the name of his first interactive theorem prover \cite{deBruijn1983}. It had a reputation for being ahead of its time, prefiguring theoretical developments that are key in modern proof assistants. But also for AutoMath proofs being unreadable (which is not unexpected in such an experimental system---especially one from a time where many popular character encodings didn't even support lowercase letters, much less any kind of symbol from ordinary mathematics).

\subsection{Summary}

We expect the transformation of mathematics set in motion by the advancement of interactive theorem proving to have far-reaching consequences. In particular
\begin{quote}
\begin{itemize}
    \item mathematics will start re-converging with the theoretical computer science;
    \item some mathematical disciplines may change language to an equivalent\footnote{In the sense of having the same scientific content, as defined by Nelson~\cite{nelson-scientific}.} one, which is more suitable for use with the new technology;
    \item new directions of research in mathematics, including those relevant to exogenous 21st century challenges such as climate change and the environment, are likely to be developed from the bottom up in new, proof-assistants-friendly ways.
\end{itemize}
\end{quote}

With this in mind, the next generation of mathematicians, regardless of their specific areas of research, will benefit from some working knowledge of proof theory, related areas of theoretical computer science, as well as good coding skills. The development and support of proof assistants, provision of help and advice regarding their use and the maintenance of proof libraries will become an important professional specialisation within mathematics. With it, the following sections of the Mathematical Subject Classification 2020 will increasingly take center stage:
\begin{quote}
\textbf{68Vxx Computer science support for mathematical research and practice}
\begin{description}
\item[68V15] Theorem proving (automated and interactive theorem provers, deduction, resolution, etc.)
\item[68V20] Formalization of mathematics in connection with theorem provers
\item[68V35] Digital mathematics libraries and repositories
\end{description}
\end{quote}

All these changes are likely to force  re-shaping of
\begin{quote}
\begin{itemize}
\item the role and place of mathematics in the society,
\item of the culture and structure of the mathematics community,
\item in publishing of mathematical research, and
\item in relations between mathematics and the state and industry.
\end{itemize}
\end{quote}

At the moment, we are not in a position to foresee any details of these future developments.

However, we have to note  the   already growing  demand for application of  \AM, or using its methods,  from the all-important safety- and -security critical industry.

However, this is easier to do in the case of mathematics education, which will be done in the next section.

Finally, we have to mention a likely impact of \AM\ on philosophy of mathematics. It is widely accepted among practicing  mathematicians, regardless of their philosophical inclinations, that the everyday functioning of mathematics is governed by a plethora of  social conventions. As it is succinctly formulated by Reuben Hersh \cite[p.29]{Persson2021}, in mathematics
\begin{quote}
[I]n practice truth is agreed on by a process of social confirmation.
\end{quote}
For most mathematicians, `truth' is equivalent to `have been proven'; but in  \AM, process of proof is no longer  a `process of social confirmation', but a technological process. What perhaps remains a social matter are possible attempts to reach a better understanding of what actually has happened inside of a computer-checked proof -- which is likely to amount to construction of another, better structured, more transparent, better annotated,  and ideally human-readable  proof. In any case, the concept of `truth' comes much closer to some form of  `absolute' than it has ever been.

Interestingly, this does not change a lot in the famous Platonism / formalism dilemma, but is likely to be hard on the social constructivism \cite{Ernest1998}.  After all, \AM\  looks as a kind of  ``virtual''  technology, but social constructivism has not been successful as a philosophy of technology  \cite{Winner1993}.

\section{Issues for mathematics education}

\subsection{The increasing role of abstractions}

Merging big parts of mathematics and computer science within  \AM\  will mean that a learner expecting to work, in the future, in \AM\  or its industrial application, will have to master a much higher level of abstract thinking.\footnote{Indeed, there are many levels of abstraction, not only in mathematics, but even more so in theoretical computer science.}  Indeed, even learning modern computing (that is, computer programming) requires a very high level of  abstraction. This view is actively promoted by Kramer \cite{Kramer2007}. A  professional computer scientist commented on Kramer's paper but preferred to stay anonymous:
\begin{quote}
``I  would caution everyone not to confuse ``mathematical thinking''  with ``the thinking done by computer scientists and programmers''.

Unfortunately, most people who are not computer scientists believe these two modes of thinking to be the same. The purposes, nature, frequency and levels of abstraction commonly used in programming are very different from those in mathematics.''
\end{quote}

 In \AM, the level of abstraction is even higher than either in  computer science, or in most areas of mathematics.

 This is a serous challenge which cuts the \AM\   away from the mainstream education of nowadays. Indeed,  in England, a modest theorem about an equivalence relation partitioning a set into equivalence classes is seen by  our university colleagues as  \textit{pons asinorum} of undergraduate abstract mathematics. Alas, many graduates from English universities obtain their bachelor’s degree in mathematics without grasping this concept.

 \subsection{The socio-economic environment}
 \label{subsec:socio-economic}

We refer the reader to papers \cite{Borovik2016,Borovik2017}; the first of them explains the changing role of mathematics in division of labour in the global economy, the second one suggests that
\begin{quote}
The current crisis in the school-level mathematics education is a sign that it reaches a bifurcation point and will inevitably split into two streams:
\begin{itemize}
\item education for a selected minority of children/young people who, in their adult lives, will be filling increasingly small share of jobs which really require
mathematics competence [for the purpose of this paper we shall call it the \textit{Deep Stream}; and
\item  awareness classes for the rest of population, end users of technology saturated by mathematics which
however will remain invisible to them (\textit{Mainstream}).
\end{itemize}
\end{quote}
In the present paper, we propose to start a systematic re-thinking of basic assumptions underpinning   mathematics education in the both streams, this will be discussed in  Sections~\ref{sec:DeepStream} and \ref{sec:Mainstream}.

Analysis of \AM\  in previous sections makes it sufficiently clear that people involved in \AM\  or in applications of \AM\  should have much better mathematical / computer science education from much earlier age, and this education should include a systematic development of  abstract thinking. This can be achieved only within the Deep Stream of school mathematics education followed by advanced level mathematical studies at universities.

However, creation of the Deep Stream immediately raises questions about the future of the rest of mathematics educations, that is, the Mainstream. These are difficult question, they will be discussed in Section~\ref{sec:Mainstream}. However, we have to emphasise that we are not talking about reforms today or tomorrow. The cycle of reproduction of mathematics
\begin{center}
 primary school -- high school -- college /  university/  -- return to school as a teacher
\end{center}
is at least 15 years long; if it continues as
\begin{center}
  university/ -- PhD studies -- postdoc -- return to university as a lecturer,
\end{center}
 it becomes at least 20 years long. Planning, preparation, and staged implementation of any serious reform of school mathematics education should be done on a comparable timescale. Education of new teachers and re-education of existing teachers, and their professional development are paramount for the success of a reform -- otherwise it will meet the same  fate as the collapsed ``New Math'' reforms of the 1960s and 1970s, with one of the more illuminating cases being Kolmogorov's reform in in Russia \cite{Borovik2021,Neretin2019}.

Some nations already have some elements of the Deep Stream education; one of them is Russia, where it is represented by specialist mathematical schools and mathematical classes, supported by a mathematical circles and mathematical competitions, etc.. We do not wish to go in the details and give only two quotes from Russian colleagues:  Sharygin \cite{Sharygin2002}
\begin{quote}
It is interesting that the Soviet system of work with mathematically gifted children, created by disinterested enthusiasts and brought, oddly enough, to the level of “know-how,” turned out to be almost the only market product of the Russian education system (not counting, of course, its final result—scientists).
\end{quote}
and Konstantinov and Semenov  \cite{Konstantinov-Semenov2021}:
\begin{quote}
[S]chools with deeper study of mathematics (mathschools) became the most important and very productive phenomenon in Russia’s education of the last decades. (p. 414)
\end{quote}
A brief introduction to the Russian tradition of specialist mathematics schools and their history can be found in \cite{Borovik2012,Gerovich2019}. It is worth mentioning that one of the leading universities in Russia, the Higher School of Economics, runs a two-year MSc programme for future  teachers of advanced level mathematics in specialist mathematics schools and classes.

Since publication of \cite{Borovik2016,Borovik2017} in 2006 and 2017, there were new surprising development which indicate that in some countries,  a split between the Deep Stream and Mainstream is becoming almost inevitable, and not because of developments within mathematics itself, but for purely socio-economic reasons---this is the theme of the next section.

\section{A systematic attack on abstract mathematical thinking}
\label{sec:Attack}

\medskip

\subsection*{This section, in its entirety, is written by A.B..  Z.K. and V.K. do not carry any responsibility for its content and views expressed.}

\medskip

\subsection{Trouble with Algebra}

As it was already discussed in Section \ref{sec:Crisis}, the cutting edge of mathematics is becoming sharper and harder than ever before. But we have to take into account a frustrating opposite development: The other end of the mathematics spectrum---the elementary mathematics in schools---is under attack from certain quarters  for being excessively abstract, too complex and remote from real life, and even worse---for being a weapon of intimidation and oppression of historically marginalised peoples.

In the USA, the mainstream mathematics education in forms proposed, say
by the Equitable Math movement \cite{Pathway2020}  or in the reform of
school mathematics in California\footnote{It is currently moving through formalities in the state legislative system, to be enacted, if approved, in July 2022.}  \cite{California2021} provide intriguing cases. Indeed, if one tries to think about the kind of socio-economic conditions these two
proposals fit the best, the answer is obvious: the system of the Universal Basic Income.

The California case deserves some special attention: this is a very technologically
developed state with world-leading hi-tech industries. So far these
industries  appear to be not much concerned about the fate of mathematics in state
schools. They are getting mathematically competent employees (and expect
to continue to get them) elsewhere. But we do not wish to preempt the
decision of the legislative assembly of California, so let us wait and see.

An Open Letter from Californian mathematicians \cite{Replace2021} against the proposed new \emph{California Math Curriculum Framework} starts with  a harsh warning:
\begin{quote}
    The proposed framework would, in effect, de-mathematize math.
\end{quote}
and adds that the Framework attempts
\begin{quote}
    ``to build a mathless Brave New World on a foundation of unsound ideology.''
\end{quote}
What is essential for the purposes of the present paper,  the new Framework, in  words of the Open Letter,
\begin{itemize}
\item ``Reject[s] ideas of natural gifts and talents'' and discourages accelerating mathematics students. 
\item Encourages keeping all students together in the same math program until the 11th grade and argues that offering differentiated programs causes student ``fragility'' and racial animosity. 
\item Rejects the longstanding goal of preparing students to take Algebra I in eighth grade, on par with high-performing foreign countries whose inhabitants will be future competitors of America's children [\dots ]—a goal explicitly part of the 1999 and 2006 Math Frameworks. 
\end{itemize}

We can add to the last point that some companies based in California \href{https://venturebeat.com/2017/02/02/15-of-facebook-employees-vulnerable-to-trumps-likely-changes-for-h-1b-visas/}{have over 15\% of their workforce imported from abroad}.

We think it is important to recognize cultural differences, and the resulting differences in the educational experiences of the students---and in full accordance with that to firmly reject the idea of intellectual barriers born of supposed ``cultural characteristics''.

If public school students are held back and not allowed to progress appropriately at their level through a variety of pathways (towards the level necessary to access their desired career paths, and to pass college admissions that put emphasis on the rigor of courses and on AP performance), private schools and tutoring services will happily come forward to fill the gap---but only for those who can afford them. Economically and socially disadvantaged students who previously had access to appropriate (standard or accelerated) courses in public middle and high school are the ones most hurt by such policies.

In one aspect of substandard orientation for waged labor, the proposed California standard imitates a failure of Canada's residential school system\footnote{While hopefully it will never come to replicate the malnutrition, disease, physical abuse, suicide, and genocide perpetrated by that tragic institution.}.
A  school inspector  wrote in 1940 in his report on a residential schools:
\begin{quote}
``Pupils who have passed through these schools possess a certain  amount  of  skills  for  which  there  is  no  demand,  and opportunities for economic self-sufficiency are neglected for want of adequate educational training.'' \cite{Stehelin1940}
\end{quote}

It is difficult to expect that the new Californian framework based on these principles  will help to educate a new generation of people able to work in mathematically-intensive industries of the future---without even mentioning future  mathematicians and computer scientists.

Ideology's invasion of   mathematics education is  a difficult issue to address, and we exclude more extreme cases of ideological attacks on mathematics from our discussion; in particular, we will not analyse \emph{Equitable Math}  \cite{Pathway2020} or make any further comments on \emph{2021 California Math Curriculum Framework} \cite{California2021}.

Instead,  we will focus on more civilised and coherently formulated arguments against ``abstract'' mathematics. We think the latter really deserve a proper discussion.

\subsection{Arguments against``abstract mathematics''}

It is interesting to see systematic attacks on teaching not only algebra, but any forms of  abstraction  in school mathematics.

An interesting example is provided by  Glen  Aikenhead \cite{Aikenhead2021,Aikenhead2021a,Aikenhead2021b} who sets

\begin{quote}
``a political agenda for liberating about a 70 percent majority of students (ages 12 to 18) from experiencing an unpleasant rite-of-passage into adulthood''. \cite{Aikenhead2021a}
\end{quote}

What is this agenda? In Aikenhead's opinion,``Platonist" (that is, abstract) mathematics is harmful for 70\% of students, ``those not interested in pursuing STEM (science, technology, engineering, mathematics) employment'' and who need instead a ``humanistic'' approach to teaching mathematics. ``Platonist" mathematics should remain only for the remaining 30\%.

What are the sins of ``Platonistic'' mathematics? First of all,

\begin{quote}
\begin{itemize}
\item[(1)] ``Platonist mathematics' controversial ontological axioms lead to deducing its mythical images of itself, maintained by privileged social power.
\item[(2)] Widespread beliefs in those myths.
\item[(3)] Socio-politico-economic power bestowed on those beliefs.
\item[(4)] Privileges gained by that social power.'' \cite{Aikenhead2021b}
\end{itemize}
\end{quote}

Aikenhead aims also
\begin{quote}
``to question the ethics of Platonist mathematics education dictating that all learners must enroll in a Platonist mathematics program imbued with its cycle of myths, when alternative mathematics programs are available.'' \cite{Aikenhead2021b}
\end{quote}
It is not clear from his papers what is the content of this "alternative mathematics", but one aspect is formulated quite clearly: this is  \textbf{not} mathematics education for future productive work in STEM areas.

Aikenhead's understanding of Platonism in mathematics is not recognisable for mathematicians: in his words, it is a kind of a religious belief that
\begin{quote}
``The universe is composed of abstract mathematical objects (AMO) that existed before humans walked the Earth.'' \cite{Aikenhead2021b}
\end{quote}
We think  most mathematicians assume no more than ``some AMOs proved to be helpful for our understanding of the (physical) universe'' and do not care  about the nature of AMOs in general.

\section{In defence of abstraction in mathematics education}
\label{sec:abstraction}

\medskip

\subsection*{This section in its entirety, is written by A.B..  Z.K. and V.K. do not carry any responsibility for its content and views expressed.}

\medskip

\subsection{Abstraction is all around us}

It suffices just to look around to realise that debates around Platonism are outdated, they have to be left where they belong: In the 19th century.

We now have a smaller, but increasingly messy, universe around us: the world of IT. And it is really built from AMOs. AMOs are hiding in each smartphone.  What is more abstract than money (ask Karl Marx), and even more so---electronic money? And there is no need even to ask this question about cryptocurrencies. As an old adage goes, money rules the world. In the modern times electronic money rules the world. Hence mathematical abstraction rules the world.

By the way, for many people (especially young) the virtual reality of social media is more important than the so-called ``real world''.

These virtual and electronic worlds of IT are becoming so complex that we need to develop some new  mathematics for understanding of how they work, what they actually do, and where they break down---and this new mathematics is \AM.

And \AM\  needs new mathematicians.

Returning to the  70\% /30\% division suggested by Aikenhead, we are afraid that it could actually happen to be much more skewed, say 95\% / 5\%.

\begin{figure}[t]
\begin{center}
(\textbf{a}) \hspace{1em}
\includegraphics[width=7.5 cm]{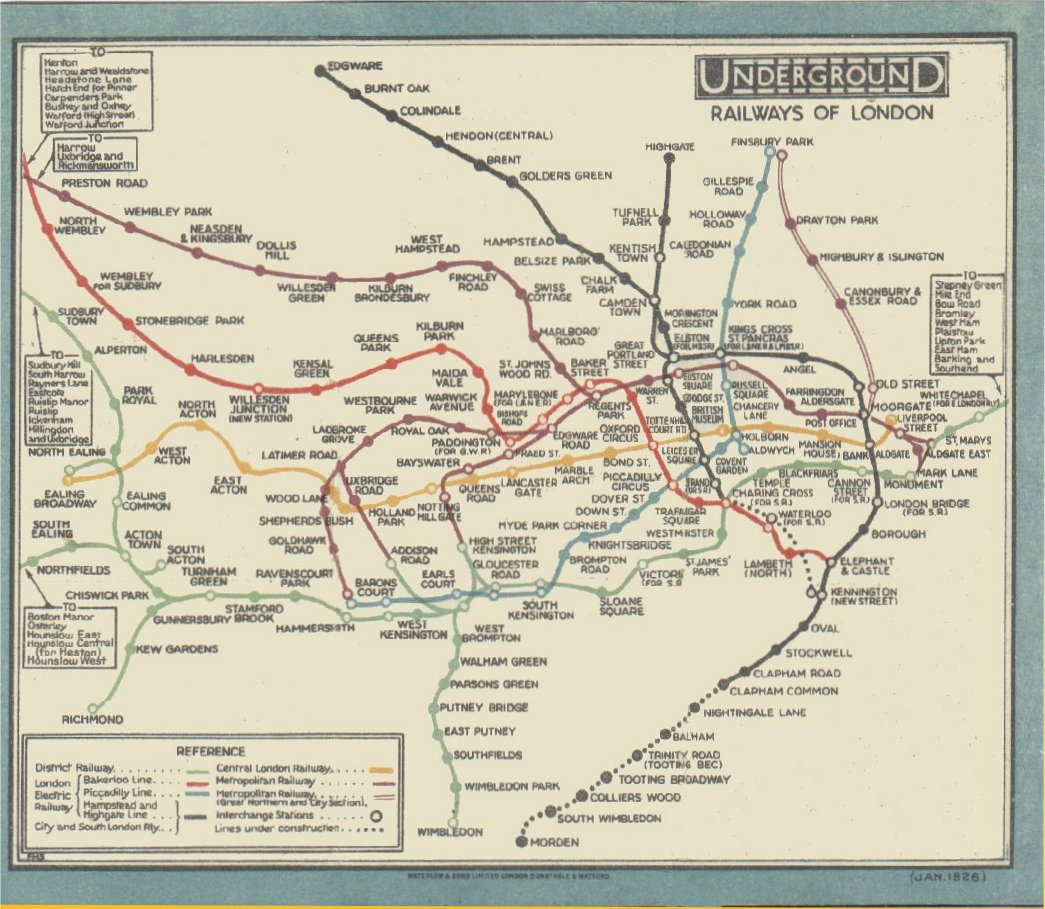}\\
\vspace{3ex}
\hspace{0.8in} (\textbf{b}) \hspace{1em}\includegraphics[width=7.5 cm]{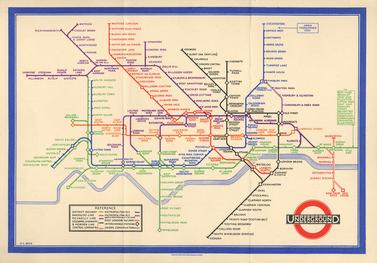}
\end{center}
\caption{Adapted from Kramer \cite{Kramer2007}.   (\textbf{a}) The London Underground map of 1928. (\textbf{b}) Harry Beck's Map of 1933, which set the template for all  further maps.  Images: \textsc{Wikimedia Commons}. \label{fig:Tube-maps}}
\end{figure}

It appears that not everyone in the mathematics education community is aware that nowadays ``abstraction'' is, first of all, a technical term of computer programming and software development, just check \href{https://en.wikipedia.org/wiki/Abstraction_(computer_science)}{Wikipedia}.
Kramer \cite{Kramer2007} gives a very useful formulation:
\begin{quote}
From the definitions of abstraction [\dots ], we focus on two particularly pertinent aspects. The first emphasises the process of removing detail to simplify and focus
attention, based on the definitions:
\begin{itemize}
\item the act of withdrawing or removing something, and
\item the act or process of leaving out of consideration one or more properties of
a complex object so as to attend to others.
\end{itemize}
The second emphasises the process of generalisation to identify the common core or
essence, based on the definitions:
\begin{itemize}
\item the process of formulating general concepts by abstracting common
properties of instances, and
\item a general concept formed by extracting common features from specific
example.
\end{itemize}
\end{quote}
Kramer then illustrates it by an example on the boundary between the information world and the real world:   maps of the London Underground,  Figure~\ref{fig:Tube-maps}(a),(b).  Map (b) became one of the iconic examples of the 20th century design. It was developed by Harry Beck from the earlier map (a) by \textit{abstraction}, removal of all unnecessary information---and this proved to be  much more efficient and convenient. As Kramer says: ``A wonderful example of the utility of abstraction''.
 Also please notice the aesthetically pleasing proportions of Beck's map, and its layout optimised for the ease of reading---it was abstraction that made it possible. And this is one of the important uses of abstraction: to present some information in a form more open to human comprehension.

\subsection{Abstraction and the real world and in mathematics}

This section uses some material from book \cite{MR2583806} and paper \cite{BOROVIK2021104410} by the first author.

We will be using the definition (or description) of mathematics from \cite{BOROVIK2021104410}:

\begin{quote}
Mathematics is the study of mental objects and constructions with reproducible properties which imitate the causality structures of the physical world, and are expressed in the human language of social interactions.
\end{quote}

The most basic elements of the causality structures of the world are schemes for expression of observations of the world so self-evident that they never mentioned in physics. For example, if you have some spoons and some forks in your cupboard and you can arrange them in pairs, with no spoon and no fork being singled out, and if you then mix spoons and forks in a box  and start matching them in pairs again, it \emph{must} be a perfect match.

Please notice the word \emph{must} -- its basic use is for expressing relations between people; please also notice that words like `must', `forces', `follows',  `defines', `holds'  etc.\ normally used for description of actions of people and relations between people, play an essential role in any mathematical narrative. What we see in the example with spoons and forks is the mathematical concept of the  one-to-one correspondence between finite sets -- as it appears ``in the wild''.  A mental construction on the top of one-to-one correspondence produces  natural numbers, arithmetic operations, and the order relation. They are interesting for their universal applicability:
\begin{quote}
 The number of my children is smaller  than the number of Galilean moons (Galilean moons are satellites of Jupiter visible from Earth via a primitive telescope or standard binoculars).
\end{quote}
This is a  statement about two groups of objects in the real world which have absolutely no ``real world'' connections between them, and which is meaningful regardless of whether it is true or false. The humble natural numbers are already  a huge abstraction.

There is another, almost folklore description of abstraction popular among computer scientists:
\begin{quote}
    Abstraction is multiple representability.
\end{quote}
Indeed, number $2$ can be represented by pairs (spoon,fork) as well as by pairs (fork, spoon) as well as by pairs of ears.

As simple as that.

\subsection{Reification}

In computer science, there is also another interesting process, \textit{reification}. From

\href{https://en.wikipedia.org/wiki/Reification_(computer_science)}{Wikipedia}:
\begin{quote}
Reification is the process by which an abstract idea about a computer program is turned into an explicit data model or other object created in a programming language. A computable/addressable object—a resource—is created in a system as a proxy for a non computable/addressable object.
\end{quote}
The results of reification are working in the IT world around us. As a technical tool, reification works in mathematics as well, although it could happen to be a more  sophisticated process than abstraction. For example, this was exactly what Sukru Yalcinkaya and the first author did in \cite{BY2018}.

Anna Sfard introduced the  term ``reification'' into mathematics education studies \cite{Sfard-Thompson1992}; she applied it to the process of objectivization of mathematical activities. The associated verb is to `reify', with the meaning ``to convert mentally into a thing, to materialize''.  The term ``reification'' is close, with some subtle differences, to that of \textit{encapsulation} \cite{Weller2004}. A detailed discussion can be found in \cite[Chapter 6]{MR2583806}. Encanpsulation is explained by an example in \cite[p. 744]{Weller2004}:
\begin{quote}
The encapsulation and de-encapsulation of process in order
to perform actions is a common experience in mathematical
thinking. For example, one might wish to add two functions
$f$ and $g$ to obtain a new function $f + g$. Thinking about doing
this requires that the two original functions and the resulting function are conceived as objects. The transformation is imagined by de-encapsulating back to the two underlying processes and coordinating them by thinking about all of the elements $x$
of the domain and all of the individual
transformations $f ( x )$ and $g ( x )$ at one time so as to obtain,
by adding, the new process, which consists of transforming
each $x$ to $f ( x ) + g ( x )$. This new process is then encapsulated
to obtain the new function $f+g$.
\end{quote}
In \cite[Chapter 6]{MR2583806} you will find also a discussion of the reverse operation, \textit{de-encapsulation}.

\section{The Deep Stream}
\label{sec:DeepStream}

\subsection{The aims of the Deep Stream}

This section develops ideas first stated in \cite{Borovik2016,Borovik2017}.

On the basis of the previous the  aims the Deep Stream could be already formulated as

\begin{quote}
developing in children skills of
\begin{itemize}
        \item identification of hidden structures in the real world, in IT, in mathematics
        \item abstraction
        \item reification and encapsulation / de-encapsulation
        \item confident handling and creation of proofs
        \item confident and fluent computer programming skills
    \end{itemize}
\end{quote}
A more detailed list of desirable outcomes, with emphasis on cognitive skills and psychological traits, is given in \cite{Borovik2017}, with a detailed discussion:
\begin{quote}
\begin{itemize}
\item[$\bullet$] ability to engage the subconscious when doing mathematics;
\item[$\bullet$] ability to share intuition;
\item[$\bullet$] ability to learn by absorption;
\item[$\bullet$] ability to compress mathematical knowledge;
\item[$\bullet$] capacity for abstract thinking;
\item[$\bullet$] being in control of their mathematics.
\end{itemize}
\end{quote}
As you can see, this sets quite challenging tasks for the Deep Stream education. However, a number of difficulties could be circumvented by a careful academic selection of students.
Khalin, Vavilov  and Yurkov \cite{Vavilov2021} share the sentiment expressed by Vladimir Rokhlin and widely accepted in Russia:
\begin{quote}
Teaching mathematics to the would-be mathematicians is infinitely easier than teaching mathematics to non-mathematicians \cite{Rokhlin2021}.
\end{quote}
For that reason we discuss here only the basic principles of the Deep Stream education.

\subsection{The Deep Stream: phase transitions}

\epigraph{We are caterpillars of angels.}{Vladimir Nabokov}

Many specialists in software and IT development working at the top of their profession passed  in their lives through 6 (six!) changes of paradigms of computer programming. Remarkably, a deep study of mathematics allows to develop skills for this kind of re-learning, changing the way of thinking:  because changes of language, of conceptual frameworks are inevitable already in learning mathematics, see \cite[Section 10]{Borovik2016} for more detail. Therefore a child in the Deep Stream needs to be gently guided through phase transitions which are similar  to  metamorphosis of a caterpillar into a butterfly.
In  \cite[Section 8]{Borovik2016}, \emph{deep mathematics education} is defined  as
\begin{quote}
Mathematics education in which every stage, starting from pre-school, is designed to fit the individual cognitive profile of the child and to serve as propaedeutics of his/her even deeper study of mathematics at later stages of education---including transition to higher level of abstraction and changes of conceptual frameworks.
\end{quote}
The concept of ``deep mathematics education'' was borrowed from Maria Droujkova, one of the leaders of American mathematics homeschooling. Her understanding of this term is, first of all, deeply human and holistic. In her own words\footnote{Private communication.},
\begin{quote}
When I use the word ``deep'' as applied to mathematics education, [\dots\ ] it means deep agency and autonomy of all participants, leading to deep personal and communal meaning and significance; as a corollary, deep individualization of every person's path; and deep psychological and technological tools to support these paths.
\end{quote}

\subsection{Where to start?}

Some suggestions are made in \cite{Borovik2017B}: From the most basic arithmetic and \textit{word problems}. That paper, in particular, reminds the principle formulated by Igor Arnold (1900--1948) \cite{Arnold1946}.
\begin{quote}
[T]eaching arithmetic involves, as a key component, the development of an \textit{ability to negotiate situations whose concrete natures represent very different relations between magnitudes and quantities}. [\dots ] The difference between the ``arithmetic'' approach to solving problems and the algebraic one is, primarily the need to make a  concrete and sensible interpretation of all the values which are used and/or which appear at any stage of the discourse.
\end{quote}
We can translate this as
\begin{quote}
The development in a child  of  an ability to see mathematical structures and relations in the world and translate them in the language of arithmetic.
\end{quote}

Krutekskii in his famous study of psychology of mathematical abilities in schoolchildren \cite{Krutetskii76} notices, as a characteristic trait of so-called ``mathematically able'' children the ability to make and use generalisations---often quite quickly. This
is one of the most basic abilities, easily detectable even at the level
of primary school: after solving a single example from a series, a
child immediately knows how to solve all examples of the same kind. We suggest that this trait can be developped in a child: for that, after solving several problems of the same type, a child should be invited to write a general solution in a toy programming language, for example, \href{https://scratch.mit.edu/}{\textsc{Scratch}}---it allows to handle integer variables. Of course, computer languages of \href{https://www.scratchjr.org/}{\textsc{ScratchJr}}  and later \href{https://scratch.mit.edu/}{\textsc{Scratch}} kind have to be introduced in  the first days of  primary school.

\subsection{Links with computer science education}

We reiterate that we have no room in this paper for development of the full curriculum for the Deep Stream. We indicate a few principal points. One of them is the need to treat mathematics and computer science as a single subject. In particular, algebra will greatly benefit from being taught and learned in parallel with a sufficiently  user-friendly programming language which includes elements of functional programming.

Many seminumerical algorithms from volume 2 of Donald Knuth's epic \emph{The art of computer programming} \cite{Knuth81} beg to be included in the course of algebra.

\section{The Mainstream mathematics education}
\label{sec:Mainstream}

\epigraph{Is there any thing whereof it may be said, See, this is new? it hath
been already of old time, which was before us.

There is no remembrance of former things; neither shall there be any
remembrance of things that are to come with those that shall come after.}{Ecclesiastes 1:10--11 KJV}

\epigraph{For mathematics, in a wilderness of tragedy and change, is a creature of the mind, born to the cry of humanity in search of an invariant reality, immutable in substance, unalterable with time.}{Cletus O. Oakley}

\epigraph{We have to do mathematics using the brain which evolved
30 000 years ago for survival in the African savanna.}{Stanislas Dehaene}

\noindent
In the situation of a spit between the Deep Stream and Mainstream educations, it is the Mainstream which represents the hardest problem. The Mainstream is much bigger and is deeply rooted in the tradition, so everyone has an opinion about it. The simplest solution is, of course, to leave it as is, do nothing. However, there is also a growing dissatisfaction with its state. For some critics, it is too weak  and insufficiently deep. We can guess that these critics could be satisfied  of the creation to the Deep Stream. However, in a growing number of countries there is also criticism of mathematics educations as being old-fashioned, teaching unnecessary stuff unrelated to ``real life'' and being excessively abstract.
A very fresh and very well founded view on this problem is formulated by Khalin, Vavilov  and Yurkov \cite{Vavilov2021} and illustrated by a fascinating textbook \cite{Vavilov2021B}. Our approach is somewhat different from theirs, although the synergy between the two papers, theirs and ours,  is remarkable.

It is a commonly held belief that  the society / nation /state expect school leavers
\begin{quote}
\begin{enumerate}
\item To be able to effectively use mathematical competence in professional activities.
\item  To be able to use basic mathematical skills and mathematical literacy to function in the world (e.g. to make \textbf{informed} economic decisions and assess the reliability of certain statements).
\item  To be able to use analytical and logical thinking to enhance lifelong learning.
\end{enumerate}
\end{quote}
The difficulty is that all parameters involved in these formulations are time-dependent. For example, the nature of professional activities changes, as well as required levels of mathematical competence -- usually they decrease, The situation with informed financial decisions at the  consumer/customer level is even worse -- banks are keen to help their customers to make these decisions via banks' apps on customers' smartphones. Analytical and logical thinking? In the pre-computer era these words were not used in mathematics education literature, their meaning is  fluid. And are people with analytical thinking really appreciated, as customers and voters, by big business and political establishment?

Despite the vagueness of these requirements, they may serve as guiding principles in reshaping the Mainstream school mathematics curriculum, but do not help to formulate selection criteria for the content of the school mathematics curriculum and  guiding principles for its   structuring.

Schweiger \cite{schweiger2006,schweiger2008} defines a \textit{fundamental mathematical idea}  as an idea that
\begin{quote}
\begin{enumerate}
 \item Appears in many fields of mathematics.
 \item Appears in mathematics on many different levels.
 \item Was important in the historical development of mathematics.
  \item Related to everyday activities.
\end{enumerate}
\end{quote}
Alas, point 4. leaves next to nothing on the list of fundamental ideas. Who extensively use, in their ``everyday life'', anything beyond the basic arithmetic (without even times tables and long division)?

Alas, the most basic questions of design of the Mainstream curriculum remain wide open.
\begin{quote}
\begin{enumerate}
 \item Which themes, which  ideas should be included in the school curriculum?

 \item To what extent and in what sequence should these concepts be studied?

 \item What are the effective ways to make mathematics curriculum connected and coherent?

 \item What is the range of real-life problems that students should be invited to discuss and solve using mathematics?

\end{enumerate}
\end{quote}
The last point deserves some special thought. A proper discussion of the most basic stuff: family budget -- could be very  painful  for people from some social strata in some countries. This might be the reason why management of personal finances is rarely discussed in the curricula in full detail.

Still, we propose for discussion  a list of questions which  appear to be unavoidable in any discussion of the  Mainstream mathematics education:

\begin{enumerate}
    \item What are the goals of the students and the goals of the education system? Answers to this question  could be different for students of different social strata and different cultural and religious backgrounds.
     \item What are the effective ways to select the content of mathematics education in order to match the goals of the students with the goals of the education system?
    \item What mathematical ideas are most fruitful in forming positive attitude to mathematics and growth mindset of the students?
    \item What are the teaching methods that minimize the chances of students to develop math anxiety?
    \item What are the ways to match the content of education and the methods of measurement and control? What ways of using technology allow to enhance the quality of grading without negative side effects?
    \item What are the routes for students who wish to transfer from the Mainstream to the Deep Stream?

\end{enumerate}

The key puzzle is the relationship between the ``basic'' and ``advanced'' cognitive skills. Is it possible to learn advanced technology without mastering the most basic concepts it is built upon? The accompanying puzzle is the relationship between training and education. Is it possible to make informed decisions without having some specialised training? Advanced digital technologies (including those with the elements of automated teaching) on the one hand provide the means for ensuring deliberation (reasonable slowness), thoroughness and the independence of the students in the learning process, on the other hand they make it possible to make a competent orchestration of the educational process with a large number of students (due to adaptive algorithms, clustering by groups and personalized learning trajectories). Is there an established methodology of teaching that allows to deal with completely different (in any reasonable sense) students inside one classroom? Is it possible to hold the attention of a large group of people, effectively synchronising their cognitive processes? Is it possible to teach without obliging anyone to know anything? The teacher is obliged to create conditions for the self-realisation of each student. Is it possible to find perfect balance between the applications and the pure ideas?

One may need to teach students to make mistakes purposefully and with pleasure. When teaching mathematics, it is necessary to develop, first of all, the ability to work with imaginary objects as with real ones. In order to enable this process, it is necessary to provide real-life processes that may serve as the working metaphors for mathematics operators. One may think of 4 different modes of thinking by analogy: the representation of physical reality by means of mathematical objects and operators, the representations of physical objects by means of other physical objects, the representations of mathematical objects by means of different mathematical objects and the representations of mathematics by means of physical analogies. It might also be helpful to ask some confusing questions: where do the concepts come from? Who was the first mathematician in human history to consider this concept? Deep understanding of mathematics is born in the process of mathematical activity, it accumulates extremely slowly and it is almost impossible to force it. When doing mathematics, it is extremely important to take into account the limitations of time and the attention span of the student.

\section{Conclusion}

This paper describes  a  tectonic shift in mathematics and mathematics education. Inevitably, it will take some time---but we believe that it is unstoppable. We have already mentioned in Section \ref{subsec:socio-economic} that a natural minimal cycle of reproduction of mathematics is 15 to 20 years long.  The  impact of the shift will be fully evident only on that scale. At this stage, it is meaningless to try to foresee  minor details, but it is it right time to start a systematic discussion of  most significant  potential developments. One of the them is perhaps inevitable split of mathematics education in two (or more) streams.

\subsection*{Authors' contributions} The authors equally contributed to the paper, with the exception of Section  \ref{sec:Proof-assistants}, written by Z.K.,  \ref{sec:Attack} and \ref{sec:abstraction}  and written  by A.B. (so Z.K. and V.K. bear no responsibility for their content), and Section \ref{sec:Mainstream}, written mostly by V.K.

\subsection*{Funding} The first author received no funding. The work of V. Kondratiev was carried out with the support of the Russian Foundation for Basic Research (Project 19-29-14217  Prospective directions and forms of the use of computer technology in school mathematics curricula).

\subsection*{Acknowledgements} A. Borovik  thanks Glen Aikenhead, Inna Capdeboscq, Gregory Cherlin, Tony Gardiner,  Konstantin Lebedev, and Alexei Muravitsky for useful discussions---but they are not responsible for any views expressed in this paper. V. Kondratiev thanks Alexei Semenov for the constant support of the study and detailed discussions.

\subsection*{Conflict of interest} The authors declare no conflict of interest. The funders had no role in the design of the study; in the writing of the manuscript, or in the decision to publish the~results.

\bibliographystyle{amsplain}
\bibliography{XXI.bib}

\providecommand{\bysame}{\leavevmode\hbox to3em{\hrulefill}\thinspace}
\providecommand{\MR}{\relax\ifhmode\unskip\space\fi MR }
\providecommand{\MRhref}[2]{%
  \href{http://www.ams.org/mathscinet-getitem?mr=#1}{#2}
}
\providecommand{\href}[2]{#2}
\begin{thebibliography}{10}

\bibitem{Pathway2020}
\emph{\href{https://equitablemath.org}{A Pathway to Equitable Math
  Instruction}}, 2020.

\bibitem{Replace2021}
Alisher~S. Abdullayev and {888 other signatories},
  \emph{\href{https://www.independent.org/news/article.asp?id=13658}{Replace
  the Proposed New California Math Curriculum Framework. Open letter to
  Governor Gavin Newsom, State Superintendent Tony Thurmond, the State Board of
  Education, and the Instructional Quality Commission}}, 2021.

\bibitem{Aikenhead2021a}
Glen~S. Aikenhead, \emph{A 21st century culture-based mathematics for the
  majority of students}, Philosophy of Mathematical Education (2021).

\bibitem{Aikenhead2021}
\bysame, \emph{Resolving conflicting subcultures within school mathematics:
  Towards a humanistic school mathematics}, Can. J. Sci. Math. Techn. Educ
  (2021).

\bibitem{Aikenhead2021b}
\bysame, \emph{School mathematics: towards ending its cycle of myths},
  Philosophy of Mathematical Education (2021).

\bibitem{Arnold1946}
Igor~Vladimirovich {Arnold}, \emph{Principles of selection and composition of
  arithmetic problems}, Izvesiya APN RSFSR \textbf{6} (1946), 8--28, In
  Russian.

\bibitem{baez-talk}
John Baez, \emph{Mathematics in the 21st century}, Topos Institute Colloquium,
  26 March 2021, 2021.

\bibitem{Borovik2012}
Alexandre Borovik, \emph{\href{https://tinyurl.com/355ac33c}{``Free maths
  schools'': some international parallels}}, The De Morgan Journal \textbf{2}
  (2012), no.~2, 23--35.

\bibitem{Borovik2021}
\bysame, \emph{The {K}olmogorov reform of mathematics education in the {USSR}},
  Submitted, 2021.

\bibitem{BOROVIK2021104410}
\bysame, \emph{A mathematician’s view of the unreasonable ineffectiveness of
  mathematics in biology}, Biosystems (2021), 104410.

\bibitem{BY2018}
Alexandre Borovik and {\c{S}}\"{u}kru Yal\c{c}{\i}nkaya, \emph{{Adjoint
  representations of black box groups} $\mathop{PSL}_2(\mathbb{F}_q)$}, Journal
  of Algebra \textbf{506} (2018), 540--591.

\bibitem{MR2583806}
Alexandre~V. Borovik, \emph{Mathematics under the {M}icroscope: {N}otes on
  {C}ognitive {A}spects of {M}athematical {P}ractice}, American Mathematical
  Society, Providence, RI, 2010. \MR{2583806 (2010k:00006)}

\bibitem{Borovik2016}
Alexandre~V. {Borovik}, \emph{Calling a spade a spade: {M}athematics in the new
  pattern of division of labour}, pp.~347--374, Springer International
  Publishing, Cham, 2016.

\bibitem{Borovik2017B}
\bysame, \emph{Economy of thought: a neglected principle of mathematics
  education}, pp.~241--265, Springer, 2017.

\bibitem{Borovik2017}
Alexandre~V. Borovik, \emph{Mathematics for makers and mathematics for users},
  pp.~309--327, Springer International Publishing, Cham, 2017.

\bibitem{CGLS2021}
Inna Capdeboscq, Daniel Gorenstein, Richard Lyons, and Ronald Solomon,
  \emph{The classification of the finite simple groups, number 9: Part v,
  chapters 1--8: Theorem $c_5$ and theorem $c_6$, stage 1}, Mathematical
  Surveys and Monographs, Amer. Math. Soc., 2021.

\bibitem{Capek1927}
Joseph {\v{C}apek} and Karel {\v{C}apek}, \emph{A {L}ong {C}at {T}ale},
  Albatros, Prague, 1927, reprinted 1996.

\bibitem{nature-scholze-2021}
Davide Castelvecchi, \emph{Mathematicians welcome computer-assisted proof in
  `grand unification' theory}, Nature \textbf{595} (2021), no.~7865, 18--19.

\bibitem{Davis-Hersh1980}
Philip Davis and Reuben Hersh, \emph{The mathematical experience},
  Birkh\"{a}ser, Boston, 1980.

\bibitem{deBruijn1983}
Nicolaas~Govert de~Bruijn, \emph{Automath, a language for mathematics},
  pp.~159--200, Springer Berlin Heidelberg, Berlin, Heidelberg, 1983.

\bibitem{Ernest1998}
Paul Ernest, \emph{Social constructivism as a philosophy of mathematics}, State
  University of New York Press, Albany, NY, 1998.

\bibitem{Gerovich2019}
Slava Gerovitch, \emph{``{W}e teach them to be free''. {S}pecialized math
  schools and the cultivation of the {S}oviet technical intelligentsia},
  Kritika: Explorations in Russian and Eurasian History \textbf{4} (2019),
  no.~20, 717--754.

\bibitem{CSFG}
Daniel Gorenstein, Richard Lyons, and Ronald Solomon, \emph{The classification
  of the finite simple groups}, Mathematical Surveys and Monographs, Amer.
  Math. Soc., 1994--2021, and continuing.

\bibitem{hales-flyspeck}
Thomas Hales, Mark Adams, Gertrud Bauer, Tat~Dat Dange, John Harrison,
  Le~Truong Hoang, Cezary Kaliszyk, Victor Magron, Sean McLaughlin, Tat~Thang
  Nguyen, Quang~Truong Nguyen, Tobias Nipkow, Steven Obua, Joseph Pleso, Jason
  Rute, Alexey Solovyev, Thi Hoai~An Ta, Nam~Trung Tran, Thi~Diep Trieu, Josef
  Urban, Ky~Vu, and Roland Zumkeller, \emph{A {formal} {proof} {of} {the}
  {Kepler} {Conjecture}}, Forum of Mathematics \textbf{5} (2017).

\bibitem{kapranov-homotopy-hypothesis}
Mikhail~M. Kapranov and Vladimir~A. Voevodsky, \emph{Infinity-groupoids and
  homotopy types}, Cahiers de Topologie et G\'eom\'etrie Diff\'erentielle
  Cat\'egoriques \textbf{32} (1991), no.~1, 29--46.

\bibitem{karagila-zf-over-nfu}
Asaf Karagila, \emph{Why use {ZF} over {NFU}?}, Mathematics Stack Exchange,
  https://math.stackexchange.com/q/193203, (version: 2012-09-09).

\bibitem{Vavilov2021}
Vladimir Khalin, Nikolai Vavilov, and Alexander Yurkov, \emph{The skies are
  falling: Mathematics for non-mathematicians}, Mathematics (2021).

\bibitem{Knuth81}
Donald~E. Knuth, \emph{The art of computer programming. {V}ol. 2, seminumerical
  algorithms.}, Addison-Wesley., 1981.

\bibitem{Kocsis2019}
Zoltan Kocsis, \emph{Development of group theory in the language of internal
  set theory}, Ph.D. thesis, Department of Mathematics, Manchester, United
  Kingdom, 2019.

\bibitem{Konstantinov-Semenov2021}
Nikolay~N. Konstantinov and Alexei~L. Semenov, \emph{Resultative education in
  mathematics schools [in russian]}, Chebyshevskii Sbornik \textbf{22} (2021),
  no.~1, 413–436.

\bibitem{Kramer2007}
Jeff Kramer, \emph{Is abstraction the key to computing?}, Commun. ACM
  \textbf{50} (2007), no.~4, 36–42.

\bibitem{Krutetskii76}
Vadim~A. Krutetskii, \emph{The psychology of mathematical abilities in
  schoolchildren}, University of Chicago Press, 1976, Translated from {R}ussian
  by J. {T}eller, edited by J.~{K}ilpatrick and I.~{W}irszup; the {R}ussian
  original published in {M}oscow in 1968.

\bibitem{nelson-scientific}
Edward Nelson, \emph{Physical reality and mathematical form}, Sankhya: The
  Indian Journal of Statistics, Series A (1961-2002) \textbf{47} (1985), no.~1,
  1--5.

\bibitem{Neretin2019}
{Yuri A.} Neretin, \emph{The kolmogorov reform of mathematical education,
  1970--1980}, \href{https://arxiv.org/pdf/1911.06108.pdf}{arXiv:1911.06108
  [math.HO]}, 2019.

\bibitem{nipkow-semantics}
Tobias Nipkow, \emph{Teaching semantics with a proof assistant: No more {LSD}
  trip proofs}, Proceedings of the 13th International Conference on
  Verification, Model Checking, and Abstract Interpretation (Berlin,
  Heidelberg), VMCAI'12, Springer-Verlag, 2012, p.~24–38.

\bibitem{Persson2021}
Ulf Persson, \emph{A conversation with reuben hersh}, Europ. Math. Soc.
  Magazine, no.~121, 20--35, But in practice truth is agreed on by a process of
  social confirmation. p. 29.

\bibitem{Rokhlin2021}
Vladimir~A. Rokhlin, \emph{Teaching mathematics to non-mathematicians}, V. A.
  Rokhlin--{M}emorial. {T}opology, geometry, and dynamics, Contemp. Math., vol.
  772, Amer. Math. Soc., Providence, RI, 2021, pp.~19--32.

\bibitem{scholze_2020}
Peter Scholze, \emph{Liquid tensor experiment}, Xena Project blog post.
  December.
  https://xenaproject.wordpress.com/2020/12/05/liquid-tensor-experiment/, Dec
  2020.

\bibitem{scholze_2021}
\bysame, \emph{Half a year of the liquid tensor experiment: Amazing
  developments}, Xena Project blog post. June.
  https://xenaproject.wordpress.com/2021/06/05/half-a-year-of-the-liquid-tensor-experiment-amazing-developments/,
  Jun 2021.

\bibitem{schweiger2006}
Fritz Schweiger, \emph{Fundamental ideas: A bridge between mathematics and
  mathematical education}, New mathematics education research and practice,
  Brill Sense, 2006, pp.~63--73.

\bibitem{schweiger2008}
\bysame, \emph{Retaining the heritage-preparing the future. fundamental ideas
  of mathematics}, ICME-11 (2008) — Monterrey (México) (Monterrey
  (México)), International Commission on Mathematical Instruction, 2008.

\bibitem{Sfard-Thompson1992}
Anna Sfard and Patrick~W. Thompson, \emph{Problems of reification:
  representations and mathematical objects}, Proceedings of the Anuual Meeting
  of the International Group for the Psychology of Mathematics
  Education–North America, Plenary Sessions (Baton Rouge LA) (D.~Kirshner,
  ed.), vol.~1, Louisiana State University, 1992, pp.~1--32.

\bibitem{Sharygin2002}
Igor Sharygin, \emph{Which ``horse'' will bring death to {R}ussian mathematics?
  [in {R}ussian]}, Otechestvennye Zapiski \textbf{3} (2002), no.~2.

\bibitem{simpson-homotopy-theory}
Carlos Simpson, \emph{Homotopy theory of higher categories: From {S}egal
  categories to n-categories and beyond}, New Mathematical Monographs,
  Cambridge University Press, 2011.

\bibitem{Stehelin1940}
E.~C. Stehelin, \emph{Inspector's report, {S}t. {B}runo {S}chool. vol. 1.
  {I}ndigenous and {N}orthern {A}ffairs {C}anada, file 777/23-5-007}, 3
  December 1940.

\bibitem{nature-szpiro}
George Szpiro, \emph{Does the proof stack up?}, Nature \textbf{424} (2003),
  no.~6944, 12--13.

\bibitem{California2021}
{T}he California Department~of Education,
  \emph{\href{https://www.cde.ca.gov/ci/ma/cf/index.asp}{2021 Revision of the
  Mathematics Framework}}, 2021.

\bibitem{lean-library}
{The mathlib Community}, \emph{The {L}ean mathematical library}, Proceedings of
  the 9th ACM SIGPLAN International Conference on Certified Programs and Proofs
  (2020).

\bibitem{schultz-wiles-mo}
Matthew Towers, \emph{Awfully sophisticated proof for simple facts},
  MathOverflow, (version: 2010-10-17).

\bibitem{Vavilov2021B}
Nikolai Vavilov, Vladimir Khalin, and Alexander Yurkov, \emph{Mathematica for a
  non-mathematitian. {E}lectronic textbook}, Moscow Centre Continuous
  Mathematics Education, Moscow, 2021, in Russian.

\bibitem{voevodsky-univalent-year}
Vladimir~A. Voevodsky, \emph{Special year on {U}nivalent {F}oundations of
  {M}athematics}, Programme Proposal, Institute for Advanced Study, 2012.

\bibitem{Weller2004}
Kirk Weller, Anne Brown, Ed~Dubinskyn, Michael McDonald, and Cynthia Stenger,
  \emph{Intimations of infinity}, Notices of AMS (2004), no.~7, 741--750.

\bibitem{Winner1993}
Langdon Winner, \emph{Social constructivism: Opening the black box and finding
  it empty}, Science as Culture \textbf{3} (1993), 427--452.

\end{thebibliography}

\end{document}